# The Variational Principle for a $\mathbb{Z}_+^N$ Action on a Hausdorff Locally Compact Space

André Caldas* Hermano Farias†

September 24, 2021

**Abstract**

We extend the definition of *topological pressure* to *locally compact Hausdorff spaces*, and we demonstrate a *variational principle* comparing the topological and measure theoretic pressures. Given a continuous $\mathbb{Z}_+^N$-action $T$ over a locally compact Hausdorff space $X$ and a continuous function vanishing at infinity $f \in C_0(X)$, we define *topological pressure* $\mathrm{P}(T, f)$ using open covers of a special type we call *admissible covers*. With this *topological pressure*, we demonstrate that

$$\mathrm{P}(T, f) = \sup_{\mu} P_\mu(T, f),$$

where the supremum is taken over all $T$-invariant probability Radon measures over $X$, and is equal to 0 when there is none.

In the last section, we present an example that illustrates why *admissible covers* are so adequate to deal with the non-compact case, while some other approaches would fail.

## 1 Introduction

In 1976, in a remarkable paper, Misiurewicz gave a short and extended proof of the so called *variational principle* [see 6]. Many attempts have been made to extend *topological pressure* to non-compact spaces. In this paper, we propose a definition for dynamical systems over locally compact Hausdorff spaces, using open covers, on the same line of the original Adler, Konheim and McAndrew idea [see 1]. Then, we successfully demonstrate the validity of the *variational principle* (theorem 4.1) for this new *topological pressure*.

Unfortunately, the authors could not find on the literature precise definitions and statements in the form and generality needed. Exactly the same way as Misiurewicz, we do not assume metrizability, and we work not with just a transform $T : X \to X$, but a $G$-action (see subsection 2.2). And we successfully drop the compactness assumption. The topological space we work with is only assumed to be Hausdorff and locally compact. This generality demanded us to redo lots of work that had already been done by others. This is specially true for many results in subsection 2.6. The bright side, however, is that we were able to put together a very didactic and self contained presentation.

*Departamento de Matemática – Universidade de Brasília, Brasil. (andrecaldas@unb.br)
†Instituto de Física – Universidade de Brasília, Brasil. (hermanofarias@unb.br)

The key concept is that of *admissible covers* (definition 2.11). We use *admissible covers* instead of *uniform structures.* Over non-compact spaces, uniform structures are not unique. Our results could be formulated in terms of the uniform structure inherited from the one-point compactification of the space we are working with. But we opted for the more intrinsic approach we present here. The idea of *admissible covers* first appeared in [3].

In section 5 we give an example that illustrates how some alternative definitions of *topological pressure* would fail to satisfy the *variational principle.* In special, defining *topological pressure* using just *finite covers* would give a value that coincides with the pressure of the compactified version of the system.

# 2 Preliminaries

In this section we recall some elementary definitions and some fundamental facts which are used in the sequel. The concept of *admissible cover* in its present form was introduced in [3]. *Admissible covers* shall be used to extend the notion of *topological pressure* to a *topological system* defined over a locally compact Hausdorff space.

## 2.1 Radon Measures

If $X$ is a topological space, denote by $C_0(X)$ the Banach space of all continuous functions $f : X \to \mathbb{R}$ *vanishing at infinity*, with the supremum norm. When $X$ is compact, this is just the space $C(X)$ of all continuous functions $f : X \to \mathbb{R}$. A *finite Radon measure* over $X$ is a finite measure such that

$$\sup \mu(K) = \mu(A) = \inf \mu(U),$$

where the supremum is taken over compact sets $K \subset A$, and the infimum, over open sets $U \supset A$. Denote by $\mathcal{M}(X)$ the normed space of finite signed radon measures over $X$, where $\|\mu\|$ is the measure's *total variation.*

Although we work with Radon measures over locally compact Hausdorff spaces, we shall need a topology over $\mathcal{M}(X)$ only when $X$ is compact. When $X$ is Hausdorff and compact, the *Riesz Representation Theorem* [4, corollary 7.18] allows the isometric identification of $\mathcal{M}(X)$ and the closed subset of positive operators in $C(X)^*$. Therefore, we can induce in $\mathcal{M}(X)$ the weak$^*$ topology from $C(X)^*$. In this topology (with a compact $X$), a net $\mu_\alpha$ converges to $\mu$ if

$$\int f \, \mathrm{d}\mu_\alpha \to \int f \, \mathrm{d}\mu$$

for every $f \in C(X)$, and the set of probability measures is compact by *Alaoglu Theorem* [4, theorem 5.18].

In general, it is not true that if $\mu_\alpha \to \mu$, then $\mu_\alpha(A) \to \mu(A)$ for every measurable set $A$. The *Portmenteau Theorem* gives a condition on the set $A$ that allows us to conclude that $\mu_\alpha(A)$ converges to $\mu(A)$. A version of the theorem for metrizable spaces and sequences of measures can be found in [2, theorem 2.1 item (v)]. We need a version for arbitrary compact Hausdorff spaces and convergent nets. All we need is *Urysohn's lemma* [4, lemma 4.32].

**Lemma 2.1.** *Let $X$ be a compact Hausdorff space and $\mu_\alpha \in \mathcal{M}(X)$ a net convergent to $\mu \in \mathcal{M}(X)$. Let $A \subset X$ be an open set, and $F \subset X$ a closed one. Then,*

$$\mu(A) \leq \liminf \mu_\alpha(A)$$
$$\limsup \mu_\alpha(F) \leq \mu(F).$$

*Proof.* Given any compact $K \subset A$, *Urysohn's lemma* implies that there exists a continuous function $f : X \to \mathbb{R}$ such that

$$I_K \leq f \leq I_A,$$

where $I_K$ and $I_A$ are the indicator functions of $K$ and $A$. Therefore,

$$\mu(K) \leq \int f \, \mathrm{d}\mu = \liminf \int f \, \mathrm{d}\mu_\alpha \leq \liminf \mu_\alpha(A).$$

By taking the supremum over $K$,

$$\mu(A) \leq \liminf \mu_\alpha(A).$$

For the closed set $F$, one can just use a similar argument, or use the fact that $F^c$ is open. □

**Proposition 2.2** (Portmanteau Theorem)**.** *Let $X$ be a compact Hausdorff space and $\mu_\alpha \in \mathcal{M}(X)$ a net convergent to $\mu \in \mathcal{M}(X)$. Let $A$ be a measurable set such that $\mu(\partial A) = 0$. Then,*

$$\mu_\alpha(A) \to \mu(A).$$

*Proof.* Notice that the hypothesis on $\partial A$ implies that

$$\mu(\mathring{A}) = \mu(A) = \mu(\overline{A}).$$

Now, the claim follows from 2.1, because

$$\begin{aligned} \mu(A) = \mu(\mathring{A}) &\leq \liminf \mu_\alpha(\mathring{A}) \\ &\leq \limsup \mu_\alpha(\overline{A}) \leq \mu(\overline{A}) = \mu(A). \end{aligned}$$

□

## 2.2 The $\mathbb{Z}_+^N$ Action

Let $\mathbb{Z}_+$ denote the set of non-negative integers. From now on, $N$ will be a fixed positive integer. Denote by $G = \mathbb{Z}_+^N$ the product semiring. The semiring $G$ will be kept fixed throughout this paper. To emphasize the similarities to the particular case where $G = \mathbb{Z}_+$, $k$, $m$, $n$, $p$ and $q$ will *always* be an element of $G$. Upper indexes indicate the components of an element of $G$. Notice that $G$ is a directed set with the relation

$$n \leq m \Leftrightarrow n^1 \leq m^1, \dots, n^N \leq m^N.$$

Thus, $G$ can be used as an index set for constructing nets. More about nets can be consulted in [7].

We shall make use of the product

$$nm = (n^1 m^1, \dots, n^N m^N).$$

Denote by

$$\Lambda(n) = \left\{ k \in G \,\middle|\, k^1 < n^1, \dots, k^N < n^N \right\},$$

and $\lambda(n) = \#\Lambda(n) = n^1 n^2 \dots n^k$. When taking limits in $n \in G$, we shall assume without further notice that no entry is zero. That is, we shall assume $n$ is "large enough" so that $\lambda(n) \neq 0$.

To subdivide $\Lambda(n)$ into smaller sets, define

$$\Lambda_n(q,k) = \left\{ p \in k + qG \,\middle|\, p + \Lambda(q) \subset \Lambda(n) \right\}$$

$$\Gamma_n(q,k) = \Lambda(n) \setminus \bigcup_{p \in \Lambda_n(q,k)} (p + \Lambda(q)).$$

In words, we can tile $\Lambda(n)$ with $q$-rectangles. The points in $\Lambda_n(q,k)$ are the bottom left corners of those tiles. While $\Gamma_n(q,k)$ are the elements of $\Lambda(n)$ not covered by those $q$-rectangles. While $\Lambda(n)$ behaves like a solid volume, $\Gamma_n(q,k)$ behaves like a thin surface. We omit $k$ when $k = 0$, and write $\Lambda_n(q)$ and $\Gamma_n(q)$, instead.

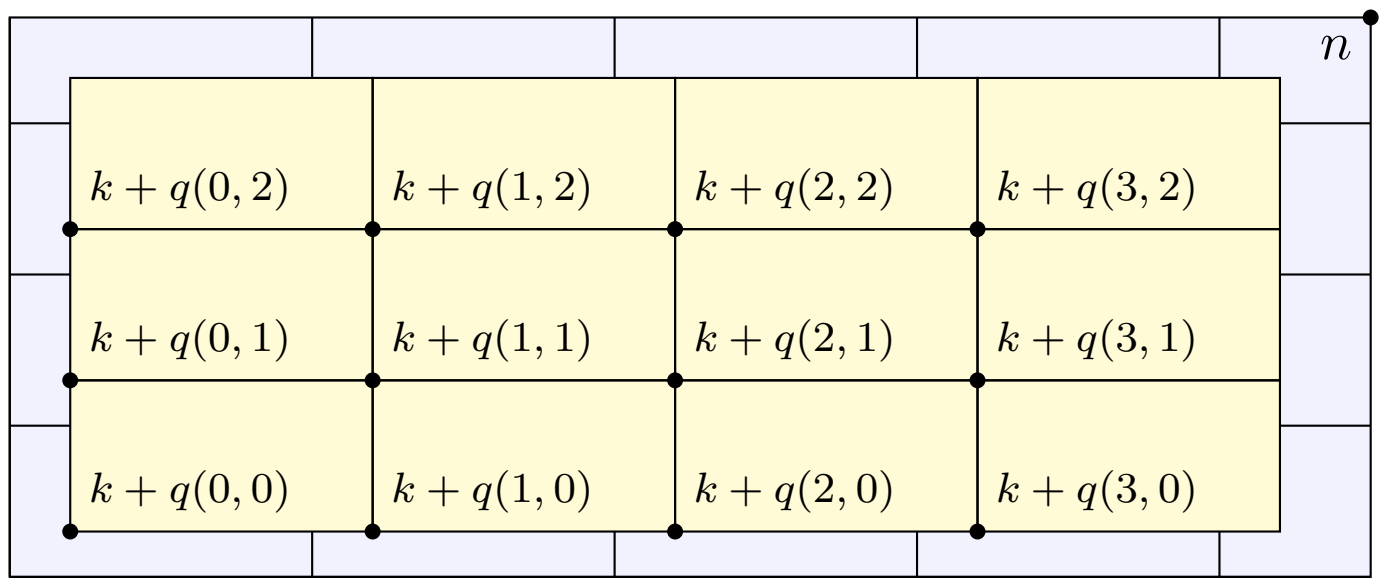

Figure 1: Partitioning the set $\Lambda(n)$. The twelve dots in the middle form $\Lambda_n(q,k)$. The twelve rectangles in the middle are of the form $p + \Lambda(q)$, for $p \in \Lambda_n(q,k)$. The outer part is $\Gamma_n(q,k)$.

The following lemma shows that when $n$ grows, the cardinality of $\Gamma_n(q,k)$ becomes small compared to that of $\Lambda(n)$.

**Lemma 2.3.** *Fixed $q, k \in G$,*

$$\frac{\#\Gamma_n(q,k)}{\lambda(n)} \to 0$$

$$\frac{\#(\Lambda(n) \,\Delta\, (q + \Lambda(n)))}{\lambda(n)} \to 0.$$

*Proof.* See figure 1. Each axis $j$ is transversal to two *"faces"* of $\Gamma_n(q,k)$. Each of these faces have *"volume"* bounded by $q^j \frac{\lambda(n)}{n^j}$. Therefore,

$$\begin{aligned}\frac{\#\Gamma_n(q,k)}{\lambda(n)} &\leq 2\sum_{j=1}^{N} q^j \frac{1}{n^j} \\ &\leq 2N\frac{\max(q)}{\min(n)} \to 0.\end{aligned}$$

In the same way,

$$\frac{\#(\Lambda(n)\,\Delta\,(q+\Lambda(n)))}{\lambda(n)} \leq 2N\frac{\max(q)}{\min(n)} \to 0.$$

□

A $G$-action $T$ over a set $X$ is a semigroup homomorphism of $G$ into the semigroup of all transformations of $X$ into itself. Instead of $T(n)$, we shall use $T^n$. That is, each

$$T^n : X \to X$$

is an application from $X$ to $X$ such that $T^{n+m} = T^n \circ T^m$, and $T^0 = \text{id}$. In spite of a $G$-action being a function

$$\begin{array}{rccc} T: & G \times X & \to & X \\ & (n,x) & \mapsto & T^n(x), \end{array}$$

we shall avoid this notation. Instead, we just write $T : X \to X$ and explicitly state that $T$ is a $G$-action. We shall also abuse the notation and write $T^{-n}$ instead of $(T^n)^{-1}$. However, when we refer to $T^n$ as a $G$-action, we mean

$$\begin{array}{rccc} T^n: & G \times X & \to & X \\ & (k,x) & \mapsto & T^{nk}(x), \end{array}.$$

For a topological (measurable) space $X$, we also say that the $G$-action is continuous (measurable) if $T^n$ is continuous (measurable) for every $n \in G$.

**Definition 2.4.** *In this work, a* topological dynamical system *is a locally compact Hausdorff topological space $X$ together with a continuous $G$-action $T : X \to X$.*

We shall omit the word *topological.* By a *dynamical system*, we mean a *topological dynamical system.* Also, since we are interested in the non-compact case, we shall assume that $X$ is not compact. This will avoid some extra care we would have to take when dealing, for example, with the *compactification* of a space that is already compact (see subsection 2.3).

**Definition 2.5.** *A* measurable dynamical system *is a measurable space $X$ together with a measurable $G$-action $T : X \to X$.*

A *topological dynamical system* is also a *measurable dynamical system* when we consider the Borel $\sigma$-algebra.

If $f : X \to \mathbb{R}$ is any function and $\Lambda \subset G$ is a finite subset of $G$, let

$$f_\Lambda = \sum_{k\in\Lambda} f \circ T^k.$$

We write $f_n$ instead of $f_{\Lambda(n)}$, for short. And if we really need to emphasize $T$, we write $f_{T,n}$ instead.

Consider the identity net $n \in G$. Let $n_\alpha$ be a universal subnet. Now, given any net $a_n$ indexed by $n \in G$, the universal subnet $n_\alpha$ induces a universal subnet $a_{n_\alpha}$. To avoid taking subnets in definitions and proofs, as well as avoid using $\limsup$ or $\liminf$, we shall fix a universal subnet $n_\alpha$. Consider, for example, the extended reals $\overline{\mathbb{R}} = \mathbb{R} \cup \{\pm\infty\}$ with its usual compact topology. If $a_n \in \mathbb{R}$, then $a_{n_\alpha}$ necessarily converges to some $a \in \bar{\mathbb{R}}$. If we show that the limit of $a_{n_\alpha}$ converges to the same $a \in \bar{\mathbb{R}}$ no matter how $n_\alpha$ was chosen, we can conclude that, in fact, $a_n \to a$. This is what we do, for example, in definition 3.3 and lemma 3.5.

## 2.3 Compactification

A compactification of a topological space $X$ is a compact Hausdorff space $Z$ containing $X$ as a dense subset, and such that $X$ has the topology induced from $Z$. A Hausdorff locally compact space $(X, \tau)$ has a *one-point compactification.* That is, take a point $\infty \notin X$ and give $X^* = X \cup \{\infty\}$ the compact Hausdorff topology given by

$$\tau^* = \tau \cup \left\{ X^* \setminus K \;\middle|\; K \subset X \text{ is compact} \right\}.$$

Since $X$ is locally compact, $X^*$ is Hausdorff. And since $X$ is Hausdorff, the topology induced in $X$ from $X^*$ is just $\tau$. In particular, [5, proposition 8.16] implies that $X^*$ has a unique uniform structure compatible with its topology. Therefore, $X$ has at least one uniform structure compatible with its topology. Then, [5, proposition 11.5] implies that $X$ is *completely regular* (see [5, definition 11.1]).

**Definition 2.6.** *With the induced topology in $X \subset Z$, a dynamical system $T : X \to X$ is a* subsystem *of $S : Z \to Z$ if $T^n(x) = S^n(x)$ for every $x \in X$. In this case, we say that $S$* extends *$T$.*

Although a locally compact Hausdorff $X$ has a one-point compactification $X^*$, it doesn't mean the dynamical system $T : X \to X$ can be continuously extended to $X^*$. Nevertheless, there does exist an extension to a (possibly much bigger) compact Hausdorff dynamical system. This is a classical result which we state and demonstrate on the next lemma, because we need a version without hypothesis like metrizability. We could just have followed the same construction as in [3, lemma 2.3]. But we opted for using the Stone-Čech compactification $\beta X$.

**Lemma 2.7.** *Suppose $T : X \to X$ is a dynamical system. Then, there exists a compactification $Z$ of $X$ such that $T$ can be continuously extended to a dynamical system $S : Z \to Z$, and such that*

$$\begin{array}{rccl} \pi: & Z & \to & X^* \\ & x & \mapsto & \pi(x) = \begin{cases} x, & x \in X \\ \infty, & x \notin X \end{cases} \end{array}$$

*is continuous.*

*Proof.* Since $X$ is *completely regular*, we can use the Stone-Čech compactification (see [9, theorem 1.11]). Let $Z = \beta X$. Notice that for each $n \in G$, we can extend the codomain of $T^n$ to $Z$. And therefore, we can continuously extend $T^n$ to $S^n : Z \to Z$. The same way, the inclusion $\iota : X \to X^*$ can be extended to $\pi : Z \to X^*$.

We need to show that for any $z \in Z \setminus X$, $\pi(z) = \infty$. But in fact, any such a $z$ is the limit of a universal net $x_\gamma \in X$. Since $x_\gamma$ does not converge to any point of $X$, it must converge to $\infty$ when seen as a net in $X^*$. Therefore,

$$\pi(z) = \lim \pi(x_\gamma) = \lim x_\gamma = \infty.$$

It only remains to show that $S$ is a $G$-action. This is a consequence of the uniqueness of the extensions $S^n$. In fact, $S^0$ is the identity, because the identity over $Z$ extends (uniquely) the identity $T^0$. Also, both $S^{n+m}$ and $S^n \circ S^m$ extend (uniquely), $T^{n+m} = T^n \circ T^m$. □

Under the conditions of lemma 2.7, $X$ is an open subset of $Z$. In fact, $X^c = \pi^{-1}(\infty)$. In this case, any Borel set of $X$ is a Borel set of $Z$. Therefore, we may restrict a Borel measure $\tilde{\mu}$ over $Z$ to the Borel sets of $X$ and produce a Borel measure over $X$, that we shall denote $\tilde{\mu}|_X$. Also, a Borel measure $\mu$ over $X$ can be extended to $Z$ by setting $\mu(X^c) = 0$.

Given a *measurable dynamical system* $T : X \to X$ and a measure $\mu \in \mathcal{M}(X)$, $\mu$ is said to be *$T$-invariant* if $\mu \circ T^{-n} = \mu$ for every $n \in G$.

**Lemma 2.8.** *Consider a measurable dynamical system $S : Z \to Z$ and a measurable subsystem $T : X \to X$. Given an $S$-invariant finite Borel measure $\mu$, then $\mu = \tilde{\mu}|_X$ is $T$-invariant.*

*Proof.* For each $n \in G$, $T^n$ is a *measurable dynamical system* in the sense of [3, lemma 2.4]. Therefore, $\mu \circ T^{-n} = \mu$. □

Instead of studying directly the system $T$, we shall look at its extension $S : Z \to Z$. However, we do not want to capture any complexity for $S$ over

$$X^c = Z \setminus X,$$

the complement of $X$ in $Z$. We say that a set $A \subset Z$ is $\pi$-measurable if $A = \pi^{-1}(B)$ for some measurable set $B \subset X^*$. A cover is $\pi$-measurable if each set in it is $\pi$-measurable.

**Definition 2.9** (One-Point Uniformly Continuous)**.** *Call a function $f : X \to \mathbb{R}$* one-point uniformly continuous*, whenever there is a continuous $\tilde{f} : X^* \to \mathbb{R}$ such that $f = \tilde{f}|_X$.*

We chose the naming *one-point uniformly continuous* because $f : X \to \mathbb{R}$ is *one-point uniformly continuous* when it is uniformly continuous with the uniform structure induced by the unique uniform structure of $X^*$. In this case, there is a unique extension $\tilde{f} : X^* \to \mathbb{R}$. And this extension induces the continuous $g = \tilde{f} \circ \pi$, from $Z$ to $\mathbb{R}$. For the dynamical systems $T : X \to X$ and $S : Z \to Z$, define

$$f_\Lambda = \sum_{k \in \Lambda} f \circ T^k$$

$$g_\Lambda = \sum_{k \in \Lambda} g \circ S^k,$$

for any finite $\Lambda \subset G$. However, $f_{\Lambda(n)} : X^* \to \mathbb{R}$ is, in principle, not defined. We use $f_n$ and $g_n$ instead of $f_{\Lambda(n)}$ and $g_{\Lambda(n)}$.

## 2.4 Coverings

A family $\mathscr{A}$ of sets whose union is equal to $X$ is a *cover* of $X$. A *finite partition* is a finite *cover* whose sets are disjoint. If $\mathscr{B} \subset \mathscr{A}$ is itself a cover of $X$, we call it a *subcover*. For $Y \subset X$, induce in $Y$ the cover

$$\mathscr{A}|_Y = \left\{ A \cap Y \;\middle|\; A \in \mathscr{A} \right\}.$$

Consider two covers of $X$, $\mathscr{A}$ and $\mathscr{B}$. If every element in $\mathscr{A}$ is contained in some $B \in \mathscr{B}$, say that $\mathscr{A}$ *refines* $\mathscr{B}$ and write $\mathscr{B} \prec \mathscr{A}$. Note that $\prec$ is a *preorder*. As usual,

$$\mathscr{A} \vee \mathscr{B} = \left\{ A \cap B \;\middle|\; A \in \mathscr{A},\, B \in \mathscr{B} \right\}$$

refines $\mathscr{A}$ and $\mathscr{B}$. We shall also use

$$\bigvee_{s \in S} \mathscr{A}_s = \left\{ \bigcap_{s \in S} A_s \;\middle|\; \forall s \in S,\, A_s \in \mathscr{A}_s \right\}$$

for a collection of covers $\mathscr{A}_s$ $(s \in S)$.

For $T : X \to X$, given a cover $\mathscr{A}$ and a subset $\Lambda \subset G$, define

$$\mathscr{A}^\Lambda = \bigvee_{k \in \Lambda} T^{-k}(\mathscr{A}).$$

We shall write $\mathscr{A}_T^\Lambda$ to indicate the corresponding dynamical system. We also write $\mathscr{A}^n$ instead of $\mathscr{A}^{\Lambda(n)}$. Notice that,

$$(\mathscr{A}_T^n)_{T^n}^m = \mathscr{A}_T^{nm}.$$

**Lemma 2.10.** *Given a subsistem $T : X \to X$ of $S : Z \to Z$. For any covering $\mathscr{Z}$ of $Z$, if we let $\mathscr{C} = \mathscr{Z}|_X$, then*

$$\mathscr{C}_T^\Lambda = \mathscr{Z}_S^\Lambda|_X.$$

*Proof.* Notice that $S$ does not take any point of $X$ into $Z \setminus X$. Therefore, for any $A \subset Z$ and $n \in G$,

$$X \cap S^{-n}(A) = X \cap S^{-n}(A \cap X) = T^{-n}(A \cap X).$$

This way,

$$\begin{aligned}
\mathscr{C}_T^\Lambda &= \bigvee_{k \in \Lambda} T^{-k}(\mathscr{C}) \\
&= \bigvee_{k \in \Lambda} (S^{-k}(\mathscr{Z}))|_X \\
&= \left( \bigvee_{k \in \Lambda} S^{-k}(\mathscr{Z}) \right)\Bigg|_X \\
&= \mathscr{Z}_S^\Lambda|_X.
\end{aligned}$$

□

## 2.5 Admissible Coverings

*Topological entropy* was introduced for compact systems in [1]. Definitions were made using open covers. For the non-compact case, using open covers without any further restriction does not work because there are covers without finite subcover. Even if one uses only finite open covers, the variational principle does not hold in general (see subsection 5.3). In [3], Caldas and Patrão successfully proposed the use of *admissible covers*.

**Definition 2.11** (Admissible Cover)**.** *In a topological space $X$, a cover $\mathscr{A}$ is said to be an* admissible cover *when it is open and at least one of its elements has compact complement. If every set has compact complement, $\mathscr{A}$ is said to be a* strongly admissible cover*.*

**Definition 2.12** (Admissible Partition)**.** *In a topological space $X$, a partition is said to be an* admissible partition *when it is finite and every element but one is compact.*

Of course, admissible partitions are measurable. The following lemma shows that *admissible partitions* are partitions that *"generate" strongly admissible covers.*

**Lemma 2.13.** *In a Hausdorff topological space $X$, if*

$$\mathscr{K} = \{K_0, \dots, K_l\}$$

*is an admissible partition where $K_1, \dots, K_l$ are all compact, then*

$$\mathscr{A} = \{K_0 \cup K_1, K_0 \cup K_2, \dots, K_0 \cup K_l\}$$

*is a strongly admissible cover.*

*Proof.* One just has to notice that $\mathscr{A}$ does cover $X$. And also, that

$$(K_0 \cup K_j)^c = \bigcup_{i \in \{1,\dots,l\} \setminus \{j\}} K_i$$

is compact for every $j = 1, \dots, l$. □

We shall need to compare the cardinality of covers and partitions.

**Lemma 2.14.** *Consider a $G$-action $T : X \to X$ over a set $X$, and a partition*

$$\mathscr{K} = \{K_0, \dots, K_l\}.$$

*Let $\mathscr{B}$ be any refinement of*

$$\mathscr{A} = \{K_0 \cup K_1, K_0 \cup K_2, \dots, K_0 \cup K_l\}.$$

*Then, for each $B \in \mathscr{B}^n$, the number of elements of $\mathscr{K}^n$ that $B$ intersects is at most $2^{\lambda(n)}$.*

*Proof.* Since $\mathscr{A} \prec \mathscr{B}$, $\mathscr{A}^n \prec \mathscr{B}^n$. Therefore, $B$ is contained in some $A \in \mathscr{A}^n$. Now,

$$A = \bigcap_{j \in \Lambda(n)} T^{-j}(K_0 \cap K_{\kappa_j}),$$

for some $\kappa \in \{1, \dots, l\}^{\Lambda(n)}$.

Therefore,

$$B \subset \bigcup_{\gamma \in \{0,1\}^{\Lambda(n)}} ( \bigcap_{k \in \Lambda(n)} T^{-k}(K_{\gamma_k \kappa_k}).)$$

Since the family $\mathscr{K}^n$ partitions $X$, $B$ does not intersect any other element of $\mathscr{K}^n$ that is not in this union. And since there is one for each $\gamma \in \{0,1\}^{\Lambda(n)}$, the claim follows. □

**Lemma 2.15.** *If $X$ is a Hausdorff locally compact topological space, then, a family $\mathscr{A}$ is an* admissible cover *of $X$ if, and only if, there is an open cover $\tilde{\mathscr{A}}$ of $X^*$, such that*

$$\mathscr{A} = \tilde{\mathscr{A}}|_X.$$

*Proof.* The restriction of an open cover of $X^*$ is admissible because any open set containing $\infty$ has compact complement contained in $X$.

On the other hand, if $\mathscr{A}$ is an *admissible cover* of $X$, then there is a set $A \in \mathscr{A}$ such that $K = X \setminus A$ is compact. Then,

$$\tilde{\mathscr{A}} = \{K^c\} \cup \mathscr{A}$$

is an open cover of $X^*$, because $X$ is open in $X^*$. □

In subsection 3.2, we shall deal with points that are $\mathscr{B}^n$-close. The points $x, y$ are $\mathscr{B}^n$-close if there is a $B \in \mathscr{B}^n$ that contains both $x$ and $y$. This is the same as saying that for all $j \in \Lambda(n)$, $T^j x$ and $T^j y$ belong to the same element $B_j \in \mathscr{B}$.

**Corollary 2.16.** *Let $X$ be a Hausdorff locally compact topological space and $f : X \to \mathbb{R}$ a one-point uniformly continuous function. Given $\varepsilon > 0$, there is an* admissible cover *$\mathscr{B}_{f,\varepsilon}$ such that whenever $x, y \in X$ are $\mathscr{B}^n_{f,\varepsilon}$-close,*

$$f_n(x) \leq f_n(y) + \lambda(n)\varepsilon.$$

*Proof.* Denote by $I_a(\varepsilon) = (a - \varepsilon/2, a + \varepsilon/2)$ the intervals of diameter $\varepsilon$ centered at $a \in \mathbb{R}$. And let $\tilde{f} : X^* \to \mathbb{R}$ be an extension of $f$. Then, $f^{-1}(I_a) = X \cap \tilde{f}^{-1}(I_a)$. So, according to lemma 2.15,

$$\mathscr{B}_{f,\varepsilon} = \left\{ f^{-1}(I_a(\varepsilon)) \,\middle|\, a \in \mathbb{R} \right\}$$

is an *admissible cover.*

If there is a $B \in \mathscr{B}^n_{f,\varepsilon}$ such that $x, y \in B$, then for $j \in \Lambda(n)$, $T^j x$ and $T^j y$ belong to a common element of $\mathscr{B}_{f,\varepsilon}$. In particular, $|f(T^j x) - f(T^j y)| < \varepsilon$. Therefore,

$$\begin{aligned} |f_n(x) - f_n(y)| &= \sum_{j \in \Lambda(n)} (|f(T^j x) - f(T^j y)|) \\ &\leq \sum_{j \in \Lambda(n)} \varepsilon = \lambda(n)\varepsilon. \end{aligned}$$

□

**Lemma 2.17.** *Let $\mathscr{A}$ be an open cover for a compact Hausdorff topological space $X$. And let $\mu$ be a Borel measure over $X$. Then, there is an open refinement $\mathscr{A} \prec \mathscr{B}$ such that*

$$B \in \mathscr{B} \Rightarrow \mu(\partial B) = 0.$$

*Proof.* For each $x \in X$, there is an $A_x \in \mathscr{A}$ with $x \in A_x$. Since $X$ is *completely regular*, there is a continuous $\varphi_x : X \to [0, 1]$ such that $\varphi_x(x) = 0$, and for any $\varepsilon \in (0, 1)$,

$$\varphi_x^{-1}([0, \varepsilon)) \subset A_x.$$

The sets $\varphi_x^{-1}(\varepsilon)$ are all disjoint and idexed by elements of the uncountable set $(0, 1)$. Since $\mu$ is finite, there is a $\varepsilon_x \in (0, 1)$ such that $\mu(\varphi_x^{-1}(\varepsilon)) = 0$. Then, the family

$$\mathscr{B} = \left\{ \varphi_x^{-1}([0, \varepsilon_x)) \;\middle|\; x \in X \right\}$$

is an open cover that refines $\mathscr{A}$. Also, the border of $\varphi_x^{-1}([0, \varepsilon_x))$ is contained in $\varphi_x^{-1}(\varepsilon_x)$, which has zero $\mu$ measure. □

## 2.6 Pressure with a Mesure

Traditionally, *pressure* has been defined only for probability measures. However, as done in [3] for entropy, extending the concept to finite measures is straightforward, and can be useful when $X$ is not compact. Although most of this section has no new results, the authors were not able to find precise statements and demonstrations in the form and generality needed for this exposition. Lemma 2.25 is a generalization of [3, lemma 2.10] and is the key result that allows us to reduce the analysis of the non-compact case to its compactification.

Didactic expositions, like [8] and [10], do not present the subject with the generality we need. While [6] does deal with $G$-actions, we could not find precise definitions and statements of needed facts by tracking the references.

We organize definitions and results in a form better suited for our purposes. For example, proposition 2.26 states explicitly that instead of general finite measurable partitions, one can simply use *admissible partitions* (definition 2.12) to calculate entropy. This, again, is not a new fact, but it is usually hidden inside the demonstration of the variational principle for compact systems.

**Definition 2.18** (Partition Entropy)**.** *Given a finite measure $\mu$ over $X$ and a finite measurable partition $\mathscr{C}$, define the* partition entropy *of $\mathscr{C}$ with respect to $\mu$ by*

$$H_\mu(\mathscr{C}) = \sum_{C \in \mathscr{C}} \mu(C) \log \frac{1}{\mu(C)}.$$

It is well known that $H_\mu(\cdot)$ is subadditive with respect to $\vee$.

**Lemma 2.19.** *Let $T : X \to X$ be a measurable dynamical system, $\mu$ a probability measure and $Q \subset G$ a finite subset. Given partitions $\mathscr{C}$ and $\mathscr{D}$ of $X$ we have that*

$$H_\mu(\mathscr{C} \vee \mathscr{D}) \leq H_\mu(\mathscr{C}) + H_\mu(\mathscr{D})$$
$$H_\mu(\mathscr{C}^Q) \leq \#Q \log \#\mathscr{C}.$$

*Proof.* The first inequality is [10, theorem 4.3, item (viii)]. Since $T$ is measurable then $T^{-q}(\mathscr{C})$ is a measurable partition. And also,

$$\#T^{-q}(\mathscr{C}) \leq \#\mathscr{C}.$$

The first inequality together with [10, corollary 4.2.1] gives

$$H_\mu(\mathscr{C}^Q) \leq \sum_{q \in Q} H_\mu(T^{-q}\mathscr{C}) \leq \sum_{q \in Q} \log \#T^{-q}(\mathscr{C}) \leq \#Q \log \#\mathscr{C}.$$

□

**Proposition 2.20.** *Given a measurable dynamical system $T : X \to X$, a $T$-invariant finite measure $\mu$ and a partition $\mathscr{C}$, the limit*

$$\lim_{n \in G} \frac{1}{\lambda(n)} H_\mu(\mathscr{C}^n)$$

*exists.*

*Proof.* Take $p \in G$. Since

$$\Lambda(n) = \Gamma_n(p) \cup \Big( \bigcup_{q \in \Lambda_n(p)} (q + \Lambda(p)) \Big),$$

we have that

$$\mathscr{C}^n = \mathscr{C}^{\Gamma_n(p)} \vee \Big( \bigvee_{q \in \Lambda_n(p)} \mathscr{C}^{q+\Lambda(p)} \Big).$$

Therefore,

$$\begin{aligned}
\frac{1}{\lambda(n)} H_\mu(\mathscr{C}^n) &\leq \frac{1}{\lambda(n)} H_\mu(\mathscr{C}^{\Gamma_n(p)}) + \frac{1}{\lambda(n)} \sum_{q \in \Lambda_n(p)} H_\mu(\mathscr{C}^{q+\Lambda(p)}) \\
&\leq \frac{\#\Gamma_n(p)}{\lambda(n)} H_\mu(\mathscr{C}) + \frac{1}{\lambda(n)} \sum_{q \in \Lambda_n(p)} H_\mu(\mathscr{C}^{\Lambda(p)}) \\
&\leq \frac{\#\Gamma_n(p)}{\lambda(n)} H_\mu(\mathscr{C}) + \frac{\#\Lambda_n(p)}{\lambda(n)} H_\mu(\mathscr{C}^p) \\
&\leq \frac{\#\Gamma_n(p)}{\lambda(n)} H_\mu(\mathscr{C}) + \frac{1}{\lambda(p)} H_\mu(\mathscr{C}^p).
\end{aligned}$$

Taking the $\limsup$ for $n$, and using lemma 2.3,

$$\limsup_{n \to \infty} \frac{1}{\lambda(n)} H_\mu(\mathscr{C}^n) \leq \frac{1}{\lambda(p)} H_\mu(\mathscr{C}^p).$$

Now, taking the $\liminf$ for $p$,

$$\limsup_{n \to \infty} \frac{1}{\lambda(n)} H_\mu(\mathscr{C}^n) \leq \liminf_{p \to \infty} \frac{1}{\lambda(p)} H_\mu(\mathscr{C}^p).$$

□

This proposition shows that the measure theoretic entropy of a dynamical system is well defined.

**Definition 2.21** (Kolmogorov-Sinai Entropy). *Given a measurable dynamical system $T : X \to X$, if $\mu$ is $T$-invariant and finite, the* partition entropy of $T$ with respect to $\mathscr{C}$ and $\mu$ *is*

$$h_\mu(T), \mathscr{C} = \lim_{n \in G} \frac{1}{\lambda(n)} H_\mu(\mathscr{C}^n).$$

*The* Kolmogorov-Sinai entropy *of $T$ with respect to $\mu$ is*

$$h_\mu(T) = \sup_{\substack{\mathscr{C}\text{: finite} \\ \text{measurable partition}}} h_\mu(T), \mathscr{C}.$$

**Definition 2.22** (Pressure). *For a measurable dynamical system $T : X \to X$ with a $T$-invariant finite measure $\mu$, and a $\mu$-integrable function $f : X \to \mathbb{R}$, the quantity*

$$P_\mu(T, f) = h_\mu(T) + \int f \, \mathrm{d}\mu$$

*is the* pressure *of $T$ with respect to $\mu$ and potential $f$. To make the notation cleaner, we also define*

$$P_\mu(T, f, \mathscr{C}) = h_\mu(T), \mathscr{C} + \int f \, \mathrm{d}\mu.$$

**Lemma 2.23.** *Given a measurable dynamical system $T : X \to X$, and a $T$-invariant finite measure $\mu$, let $f : X \to \mathbb{R}$ be an integrable function and $c \in \mathbb{R}$, a real number. Then,*

$$P_\mu(T, f + c) = P_\mu(T, f) + c\mu(X).$$

*Proof.* This is an obvious consequence of $\int (f + c) \, \mathrm{d}\mu = \int f \, \mathrm{d}\mu + c\mu(X)$. □

Usually, *entropy* and *pressure* are defined using *probability measures*. In the non-compact case, however, the weak-$*$ limit of a probability measure might not have total measure equal to 1. The following lemma helps us to deal with this fact.

**Lemma 2.24.** *Let $T : X \to X$ be a measurable dynamical system, $\mu$ a finite $T$-invariant measure and $f : X \to \mathbb{R}$ a $\mu$-integrable function. Then, for $\alpha \geq 0$ and any measurable finite partition $\mathscr{C}$,*

$$P_{\alpha\mu}(T, f, \mathscr{C}) = \alpha P_\mu(T, f, \mathscr{C}).$$

*In special,*

$$P_{\alpha\mu}(T, f) = \alpha P_\mu(T, f).$$

*Proof.* The integral part is evident, since

$$\int f \, \mathrm{d}(\alpha\mu) = \alpha \int f \, \mathrm{d}\mu.$$

The entropy part is basically [3, lemma 2.9], but we shall demonstrate it here because the setting is more general. For any finite partition $\mathscr{C}$,

$$\begin{aligned} \frac{1}{\lambda(n)} \sum_{C \in \mathscr{C}^n} \alpha\mu(C) \log \frac{1}{\alpha\mu(C)} &= \frac{\alpha}{\lambda(n)} \sum_{C \in \mathscr{C}^n} \mu(C)(\log \frac{1}{\alpha} + \log \frac{1}{\mu(C)}) \\ &= \frac{\alpha}{\lambda(n)} \mu(X) \log \frac{1}{\alpha} + \frac{\alpha}{\lambda(n)} \sum_{C \in \mathscr{C}^n} \mu(C) \log \frac{1}{\mu(C)}. \end{aligned}$$

Now, just take the limit on $n$. □

If the system $T : X \to X$ can be extended to $S : Z \to Z$, we can use $S$ to help us calculate the pressure of $T$. In the next lemma, we calculate the entropy of $S$ just for partitions that are $\pi$-measurable. That is, partitions composed of sets $\pi^{-1}(B)$, where $B \subset X^*$ is measurable.

**Lemma 2.25.** *Let $T : X \to X$ be a dynamical system and $S : Z \to Z$ the compact extension described in lemma 2.7. If $\tilde{\mu}$ is a finite $S$-invariant measure, then $\mu = \tilde{\mu}|_X$ is $T$-invariant and*

$$h_\mu(T) = \sup_{\substack{\mathscr{Z}:\ \pi\text{-measurable} \\ \text{finite partition}}} h_{\tilde{\mu}}(S), \mathscr{Z}.$$

*If $f \in C_0(X)$, then*

$$P_\mu(T, f) = \sup_{\substack{\mathscr{Z}:\ \pi\text{-measurable} \\ \text{finite partition}}} P_{\tilde{\mu}}(S, f \circ \pi, \mathscr{Z}).$$

*Proof.* The $T$-invariance of $\mu$ over $X$ is just lemma 2.8. Now, if we use lemmas 2.8 and 2.10, and divide by $\lambda(n)$ instead of dividing by $n$, [3, lemma 2.10] works verbatim to show that

$$h_\mu(T) = \sup_{\substack{\mathscr{Z}:\ \pi\text{-measurable} \\ \text{finite partition}}} h_{\tilde{\mu}}(S), \mathscr{Z}.$$

The final inequality follows from

$$\int f \circ \pi \, \mathrm{d}\tilde{\mu} = \int f \, \mathrm{d}\mu + f(\infty)\tilde{\mu}(\{\infty\}) = \int f \, \mathrm{d}\mu.$$

□

To calculate the pressure, we can restrict ourselves to *admissible partitions* (definition 2.12). This fact is usualy not explicitly stated as we do here, in proposition 2.26. But it is not a new fact, as it is usually embedded inside the demonstration of the *variational principle* for the compact case. For example, see [8] or [10].

**Proposition 2.26.** *Given a measurable dynamical system $T : X \to X$, $T$-invariant finite Radon measure $\mu$ and a $\mu$-integrable function $f : X \to \mathbb{R}$,*

$$P_\mu(T, f) = \sup_{\substack{\mathscr{K}:\ \text{admissible} \\ \text{partition}}} P_\mu(T, f, \mathscr{K}).$$

For the proof of proposition 2.26, we need the concept of *conditional entropy.* The proof will be presented after some preparation.

**Definition 2.27** (Conditional Entropy)**.** *If $\mu$ is a probability measure and $\mathscr{C}$ and $\mathscr{D}$ are two finite measurable partitions, the* conditional entropy *is the expected value*

$$H_\mu(\mathscr{C} \mid \mathscr{D}) = \sum_{D \in \mathscr{D}} \mu(D) H_{\mu(\cdot \mid D)}(\mathscr{C}).$$

Conditional entropy will be used to replace a measurable finite partition by an *admissible partition*. Lemma 2.28 gives a boundary for how much the entropy of a partition $\mathscr{K}$ might underestimate the entropy for another partition $\mathscr{C}$. And proposition 2.29 shows that the estimated difference can be made arbitrarily small if $\mathscr{K}$ is properly chosen.

**Lemma 2.28.** *Given a measurable dynamical system $T : X \to X$ with $T$-invariant probability measure $\mu$, then,*

$$h_\mu(T), \mathscr{C} \leq h_\mu(T), \mathscr{K} + H_\mu\big(\mathscr{C} \,\big|\, \mathscr{K}\big)$$

*for any two measurable finite partitions $\mathscr{C}$ and $\mathscr{K}$.*

*Proof.* We basically repeat the proof of [10, theorem 4.12, item (iv)]. Using [10, theorem 4.3, items (vii), (v) and (ix)],

$$\begin{aligned}
\frac{1}{\lambda(n)} H_\mu(\mathscr{C}^n \vee \mathscr{K}^n) &= H_\mu\left( \bigvee_{k \in \Lambda(n)} T^{-k}(\mathscr{C}) \,\middle|\, \bigvee_{q \in \Lambda(n)} T^{-q}(\mathscr{K}) \right) \\
&\leq \sum_{k \in \Lambda(n)} H_\mu\left( T^{-k}(\mathscr{C}) \,\middle|\, \bigvee_{q \in \Lambda(n)} T^{-q}(\mathscr{K}) \right) \\
&\leq \sum_{k \in \Lambda(n)} H_\mu\Big( T^{-k}(\mathscr{C}) \,\Big|\, T^{-k}(\mathscr{K}) \Big) \\
&= \lambda(n) H_\mu\big(\mathscr{C} \,\big|\, \mathscr{K}\big).
\end{aligned}$$

And [10, theorem 4.3, items (iv) and then (ii)] imply that

$$\begin{aligned}
\frac{1}{\lambda(n)} H_\mu(\mathscr{C}^n) &\leq \frac{1}{\lambda(n)} H_\mu(\mathscr{C}^n \vee \mathscr{K}^n) \\
&= \frac{1}{\lambda(n)} \left( H_\mu(\mathscr{K}^n) + H_\mu\Big(\mathscr{C}^n \,\Big|\, \mathscr{K}^n\Big) \right) \\
&\leq \frac{1}{\lambda(n)} H_\mu(\mathscr{C}^n) + H_\mu\big(\mathscr{C} \,\big|\, \mathscr{K}\big)
\end{aligned}$$

from where the result follows. ∎

Now, we need to show that we can choose an *admissible partition* $\mathscr{K}$, such that $H_\mu\big(\mathscr{C} \,\big|\, \mathscr{K}\big)$ is as small as one might require. In the following proposition, Radon measures' regularity will allow us to build up such an *admissible partition.*

**Proposition 2.29.** *Let $\mathscr{C} = \{C_1, \dots, C_l\}$ be a measurable partition. If $\mathscr{K} = \{K_0, K_1, \dots, K_l\}$ is such that $K_j \subset C_j$ for every $j = 1, \dots, l$, then*

$$h_\mu(T), \mathscr{C} \leq h_\mu(T), \mathscr{K} + \mu(K_0) \log l.$$

*Proof.* According to lemma 2.28, we just have to show that

$$H_\mu\big(\mathscr{C} \,\big|\, \mathscr{K}\big). \leq \mu(K_0) \log l.$$

Just notice that for every $C \in \mathscr{C}$ and $j \neq 0$, $\mu\big(C \,\big|\, K_j\big)$ is either 0 or 1. Therefore, for $j = 1, \dots, l$,

$$H_{\mu(\cdot \,|\, K_j)}(\mathscr{C}) = 0.$$

And [8, lemma 9.1.3] implies that $H_{\mu(\cdot\,|\,K_0)}(\mathscr{C}) \leq \log l$. Now, apply definition 2.27. ☐

We are now, ready do demonstrate proposition 2.26.

*Proof (proposition 2.26).* Using lemma 2.24, we can assume $\mu$ is a probability. An *admissible partition* $\mathscr{K}$, besides being finite by definition, it is also measurable, because compact sets are measurable. Therefore,

$$\sup_{\substack{\mathscr{K}:\text{ admissible}\\\text{partition}}} h_\mu(\mathscr{K}), T \leq \sup_{\substack{\mathscr{C}:\text{ finite}\\\text{measurable partition}}} h_\mu(\mathscr{C}), T = h_\mu(T).$$

To complete the demonstration, given any measurable finite partition $\mathscr{C} = \{C_1, \dots, C_l\}$, we find an *admissible partition* $\mathscr{K}$ satisfying proposition 2.29 and such that $\mu(K_0)$ is arbitrarily small. To this end, given $\varepsilon > 0$, since $\mu$ is Radon, choose a compact $K_j \subset C_j$ for each $j = 1, \dots, l$, such that

$$\mu(C_j \setminus K_j) \leq \frac{\varepsilon}{l \log l}.$$

Since $K_0 = (K_1 \cup \cdots \cup K_l)^c$,

$$\mu(K_0) \log l = \sum_{j=1}^{l} \mu(C_j \setminus K_j) \leq \varepsilon.$$

Now, proposition 2.29 implies that

$$\varepsilon + \sup_{\substack{\mathscr{K}:\text{ admissible}\\\text{partition}}} h_\mu(\mathscr{K}), T \geq \sup_{\substack{\mathscr{C}:\text{ finite}\\\text{measurable partition}}} h_\mu(\mathscr{C}), T = h_\mu(T).$$

Since $\varepsilon > 0$ was arbitrary, the result follows. ☐

A very well known result is the exponential-like behaviour of the pressure.

**Lemma 2.30.** *Let $T : X \to X$ be a measurable $G$-action. Given a $T$-invariant measure $\mu$ and a $\mu$-integrable function $f : X \to \mathbb{R}$, and given $m \in G$,*

$$P_\mu(T^m, f) = \lambda(m) P_\mu(T, f).$$

*Proof.* We notice that $(\mathscr{C}_T^m)_{T^m}^n = \mathscr{C}_T^{nm}$ for any partition $\mathscr{C}$. Now,

$$\frac{\lambda(m)}{\lambda(nm)} \log H_\mu(\mathscr{C}_T^{nm}) = \frac{1}{\lambda(n)} \log H_\mu((\mathscr{C}_T^m)_{T^m}^n).$$

Taking the limit for $n$,

$$\lambda(m) h_\mu(T), \mathscr{C} = h_\mu(T^m), \mathscr{C}_T^m.$$

Then,

$$\begin{aligned}
\lambda(m) h_\mu(T) &= \lambda(m) \sup_{\mathscr{C}} h_\mu(T), \mathscr{C} \\
&= \sup_{\mathscr{C}} h_\mu(T^m), \mathscr{C}_T^m \\
&= \sup_{\mathscr{D}} h_\mu(T^m), \mathscr{D} \\
&= h_\mu(T^m),
\end{aligned}$$

where the third equation follows because $\mathscr{C} \prec \mathscr{C}_T^m$.

To finalize the proof, notice that since $\mu$ is $T$-invariant,

$$\int f_m \, \mathrm{d}\mu = \lambda(m) \int f \, \mathrm{d}\mu.$$

□

Next, we present an upper bound for calculating the *pressure* that also motivates the definition of *topological pressure*. First, notice that

$$P_\mu(T, f) = \sup_{\mathscr{C}} \lim_{n \in G} \left( \int f \, \mathrm{d}\mu + \frac{1}{\lambda(n)} H_\mu(\mathscr{C}^n) \right),$$

where the supremum is taken over every measurable finite partition $\mathscr{C}$. Again, it is not a new result. It is basically [6, equation (6)]. For the case $N = 1$, it is done, for example, in [8, section 10.4.1].

**Lemma 2.31.** *Let $T : X \to X$ be a measurable dynamical system, $\mu$ a $T$-invariant probability measure, and $f : X \to \mathbb{R}$ an integrable function. Then, for every finite measurable partition $\mathscr{C}$,*

$$\int f \, \mathrm{d}\mu + \frac{1}{\lambda(n)} H_\mu(\mathscr{C}^n) \leq \frac{1}{\lambda(n)} \log \sum_{C \in \mathscr{C}^n} \sup \mathrm{e}^{f_n(C)}.$$

*Proof.* From the $T$-invariance of $\mu$,

$$\int f \, \mathrm{d}\mu = \frac{1}{\lambda(n)} \int f_n \, \mathrm{d}\mu.$$

Therefore,

$$\begin{aligned}
\int f \, \mathrm{d}\mu + \frac{1}{\lambda(n)} H_\mu(\mathscr{C}^n) &= \frac{1}{\lambda(n)} \left( \int f_n \, \mathrm{d}\mu + \sum_{C \in \mathscr{C}^n} \mu(C) \log \frac{1}{\mu(C)} \right) \\
&\leq \frac{1}{\lambda(n)} \sum_{C \in \mathscr{C}^n} \left( \mu(C) \sup f_n(C) + \mu(C) \log \frac{1}{\mu(C)} \right) \\
&= \frac{1}{\lambda(n)} \log \exp \left( \sum_{C \in \mathscr{C}^n} \mu(C) \left( \sup f_n(C) + \log \frac{1}{\mu(C)} \right) \right) \\
&\leq \frac{1}{\lambda(n)} \log \sum_{C \in \mathscr{C}^n} \mu(C) \exp \left( \sup f_n(C) + \log \frac{1}{\mu(C)} \right) \\
&\leq \frac{1}{\lambda(n)} \log \sum_{C \in \mathscr{C}^n} \mathrm{e}^{\sup f_n(C)}.
\end{aligned}$$

□

# 3 Topological Pressure

Let us now generalize and investigate the topological pressure. Many of the original definitions and results are not related to the compactness of the topological space over which the dynamical systems is defined. The central instrument will be that of *admissible covers*.

## 3.1 Using Admissible Covers

Lemma 2.31 suggests a definition for topological pressure.

**Definition 3.1.** *Given $f : X \to \mathbb{R}$ and a cover $\mathscr{A}$ of a set $X$, define*

$$\mathrm{Q}_n(T, f, \mathscr{A}) = \inf \left\{ \sum_{A \in \mathscr{A}'} \inf \mathrm{e}^{f_n(A)} \;\middle|\; \mathscr{A}' \text{ is a subcover of } \mathscr{A}^n \right\}$$

$$\mathrm{P}_n(T, f, \mathscr{A}) = \inf \left\{ \sum_{A \in \mathscr{A}'} \sup \mathrm{e}^{f_n(A)} \;\middle|\; \mathscr{A}' \text{ is a subcover of } \mathscr{A}^n \right\}.$$

The role played by $\mathrm{Q}_n(T, f, \mathscr{A})$ and $\mathrm{P}_n(T, f, \mathscr{A})$ in definition 3.1 is analogous to that of $N(\mathscr{A}^n)$ when we define topological entropy (see [3]). In fact,

$$\mathrm{Q}_n(T, 0, \mathscr{A}) = \mathrm{P}_n(T, 0, \mathscr{A}) = N(\mathscr{A}^n).$$

The following lemma shows that $\mathrm{Q}_n(T, f, \mathscr{A})$ has a property very simmilar to that of $N(\mathscr{A}^n)$.

**Lemma 3.2.** *If $\mathscr{A} \prec \mathscr{B}$, then, for any $f : X \to \mathbb{R}$, and any $n \in G$,*

$$\mathrm{Q}_n(T, f, \mathscr{A}) \leq \mathrm{Q}_n(T, f, \mathscr{B}).$$

*Proof.* Notice that $\mathscr{A} \prec \mathscr{B}$ implies $\mathscr{A}^n \prec \mathscr{B}^n$.

For every $B \in \mathscr{B}^n$, there is an $A_B \in \mathscr{A}^n$ such that $B \subset A_B$. In this case,

$$\inf \mathrm{e}^{f(A_B)} \leq \inf \mathrm{e}^{f(B)}.$$

Notice that for every subcover $\mathscr{B}'$ of $\mathscr{B}^n$,

$$\mathscr{A}' = \left\{ A_B \;\middle|\; B \in \mathscr{B}' \right\}$$

is a subcover of $\mathscr{A}^n$. Therefore,

$$\begin{aligned} \mathrm{Q}_n(T, f, \mathscr{A}) &\leq \sum_{A \in \mathscr{A}'} \inf \mathrm{e}^{f(A)} \\ &\leq \sum_{B \in \mathscr{B}'} \inf \mathrm{e}^{f(A_B)} \\ &\leq \sum_{B \in \mathscr{B}'} \inf \mathrm{e}^{f(B)}. \end{aligned}$$

The result follows if we take the infimum over every subcover $\mathscr{B}'$ of $\mathscr{B}^n$. □

We use the universal subnet $n_\alpha \in G$ we have fixed in subsection 2.2 to define pressure (definition 3.6). Later on, propositions 3.4 and 3.11 demonstrate that the definition actually does not depend on the chosen subnet $n_\alpha$.

**Definition 3.3.** *Given $f : X \to \mathbb{R}$ and a cover $\mathscr{A}$ of a set $X$, define*

$$\mathrm{P}(T, f, \mathscr{A}) = \lim \frac{1}{\lambda(n_\alpha)} \log \mathrm{P}_{n_\alpha}(T, f, \mathscr{A}).$$

$$\mathrm{Q}(T, f, \mathscr{A}) = \lim \frac{1}{\lambda(n_\alpha)} \log \mathrm{Q}_{n_\alpha}(T, f, \mathscr{A})$$

As in [3], we shall restrict our attention to *admissible covers.*

**Proposition 3.4.** *Given a $G$-action $T : X \to X$, a bounded function $f : X \to \mathbb{R}$ and a cover $\mathscr{A}$ of $X$, the limit*

$$\lim_{n\in G} \frac{1}{\lambda(n)} \log \mathrm{P}_n(T, f, \mathscr{A})$$

*exists.*

*Proof.* Take $p \in G$ big enough such that $\lambda(p) \neq 0$. Let $q \in G$ be the greatest element of $G$ such that $pq \leq n$. In this case,

$$\lambda(q) = \#\Lambda_n(p) \quad \text{and} \quad \frac{\lambda(q)}{\lambda(n)} \leq \frac{1}{\lambda(p)}.$$

Since

$$\Lambda(n) = \Gamma_n(p) \cup \Big( \bigcup_{q\in\Lambda_n(p)} (q + \Lambda(p))\Big)$$

is a disjoint union, it follows that

$$\begin{aligned}\mathscr{A}^n &= \Big( \bigvee_{q\in\Lambda_n(p)} (\mathscr{A}^{q+\Lambda(p)})\Big) \vee \mathscr{A}^{\Gamma_n(p)} \\ &= \Big( \bigvee_{q\in\Lambda_n(p)} T^{-q}(\mathscr{A}^p)\Big) \vee \mathscr{A}^{\Gamma_n(p)}.\end{aligned}$$

Also,

$$\begin{aligned}f_n &= \Big( \sum_{q\in\Lambda_n(p)} f_{q+\Lambda(p)}\Big) + f_{\Gamma_n(p)} \\ &= \Big( \sum_{q\in\Lambda_n(p)} f_p \circ T^q\Big) + f_{\Gamma_n(p)}.\end{aligned}$$

So, given a subcover $\mathscr{B} \subset \mathscr{A}^p$, let

$$\mathscr{D} = \Big( \bigvee_{q\in\Lambda_n(p)} T^{-q}\mathscr{B}\Big) \vee \mathscr{A}^{\Gamma_n(p)}.$$

Then, $\mathscr{D}$ is a subcover of $\mathscr{A}^n$. And any $D \in \mathscr{D}$ is of the form

$$D = \Big( \bigcap_{q\in\Lambda_n(p)} T^{-q}B_{q,D}\Big) \cap C_D,$$

for $B_{q,D} \in \mathscr{B}$ and $C_D \in \mathscr{A}^{\Gamma_n(p)}$. Therefore,

$$\begin{aligned}
\mathrm{P}_n(T,f,\mathscr{A}) &\leq \sum_{D\in\mathscr{D}} \sup \mathrm{e}^{f_n(D)} \\
&= \sum_{D\in\mathscr{D}} \sup(\mathrm{e}^{f_{\Gamma_n(p)}(C_D)} \prod_{q\in\Lambda_n(p)} \mathrm{e}^{f_p(T^q(T^{-q}B_{q,D}))}) \\
&\leq \sum_{D\in\mathscr{D}} ((\sup \mathrm{e}^{f_{\Gamma_n(p)}(C_D)}) \prod_{q\in\Lambda_n(p)} \sup \mathrm{e}^{f_p(B_{q,D})}) \\
&\leq \mathrm{e}^{\#\Gamma_n(p)\|f\|_\infty} \sum_{D\in\mathscr{D}} ( \prod_{q\in\Lambda_n(p)} \sup \mathrm{e}^{f_p(B_{q,D})}) \\
&\leq \mathrm{e}^{\#\Gamma_n(p)\|f\|_\infty} \prod_{q\in\Lambda_n(p)} ( \sum_{D\in\mathscr{D}} \sup \mathrm{e}^{f_p(B_{q,D})}) \\
&\leq \mathrm{e}^{\#\Gamma_n(p)\|f\|_\infty} \prod_{q\in\Lambda_n(p)} ( \sum_{B\in\mathscr{B}} \sup \mathrm{e}^{f_p(B)}) \\
&= \mathrm{e}^{\#\Gamma_n(p)\|f\|_\infty} ( \sum_{B\in\mathscr{B}} \sup \mathrm{e}^{f_p(B)})^{\#\Lambda_n(p)}
\end{aligned}$$

Since $\mathscr{B}$ is an arbitrary subcover of $\mathscr{A}^p$,

$$\mathrm{P}_n(T,f,\mathscr{A}) \leq \mathrm{e}^{\#\Gamma_n(p)\|f\|_\infty}\, \mathrm{P}_p(T,f,\mathscr{A})^{\lambda(q)}.$$

So,

$$\begin{aligned}
\frac{1}{\lambda(n)} \log \mathrm{P}_n(T,f,\mathscr{A}) &\leq \frac{\#\Gamma_n(p)}{\lambda(n)}\|f\|_\infty + \frac{\lambda(q)}{\lambda(n)} \log \mathrm{P}_p(T,f,\mathscr{A}) \\
&\leq \frac{\#\Gamma_n(p)}{\lambda(n)}\|f\|_\infty + \frac{1}{\lambda(p)} \log \mathrm{P}_p(T,f,\mathscr{A}).
\end{aligned}$$

Taking the lim sup for $n$ and using lemma 2.3,

$$\limsup_{n\in G} \frac{1}{\lambda(n)} \log \mathrm{P}_n(T,f,\mathscr{A}) \leq \frac{1}{\lambda(p)} \log \mathrm{P}_p(T,f,\mathscr{A}).$$

Now, taking the lim inf for $p$ gives the existence of the limit. □

The family of all admissible covers is preordered by $\prec$. Therefore, we can talk about the existence or not of the limits

$$\lim_{\mathscr{A}\text{: admissible cover}} \mathrm{P}(T,f,\mathscr{A}) \qquad \text{and} \qquad \lim_{\mathscr{A}\text{: admissible cover}} \mathrm{Q}(T,f,\mathscr{A}).$$

A consequence of lemma 3.2 is that this later limit exists and equals the supremum of $\mathrm{Q}(T,f,\mathscr{A})$ over all *admissible covers*. The case of $\mathrm{P}(T,f,\mathscr{A})$ is a little more delicate.

**Lemma 3.5.** *For a $G$-action $T : X \to X$, and a one-point uniformly continuous function $f : X \to \mathbb{R}$,*

$$\begin{aligned}
\lim_{\mathscr{A}\text{: admissible cover}} \mathrm{P}(T,f,\mathscr{A}) &= \lim_{\mathscr{A}\text{: admissible cover}} \mathrm{Q}(T,f,\mathscr{A}) \\
&= \sup_{\mathscr{A}\text{: admissible cover}} \mathrm{Q}(T,f,\mathscr{A}).
\end{aligned}$$

*Proof.* Let $\mathscr{B} = \mathscr{B}_{f,\varepsilon}$ from corollary 2.16. According to the corollary, $\mathscr{B}$ is an *admissible cover*. Notice that for any $B \in \mathscr{B}^n$,

$$\sup f(T^j B) \leq \inf f(T^j B) + \varepsilon \quad (j \in \Lambda(n)).$$

And this means that

$$\sup f_n(B) \leq \inf f_n(B) + \lambda(n)\varepsilon.$$

Since we are taking limits, we might consider only admissible covers $\mathscr{A}$ such that $\mathscr{B} \prec \mathscr{A}$. We shall demonstrate that

$$\mathrm{Q}(T, f, \mathscr{A}) \leq \mathrm{P}(T, f, \mathscr{A}) \leq \mathrm{Q}(T, f, \mathscr{A}) + \varepsilon. \tag{1}$$

Since the limit exists for $\mathrm{Q}(T, f, \mathscr{A})$ and equals $\sup_{\mathscr{C}} \mathrm{Q}(T, f, \mathscr{C})$; and since $\varepsilon > 0$ is arbitrary, the result shall follow from this inequality.

Any $A \in \mathscr{A}^n$ is contained in some $B \in \mathscr{B}^n$. Therefore,

$$\sup f_n(A) \leq \inf f_n(A) + \lambda(n)\varepsilon.$$

Therefore,

$$\mathrm{Q}_n(T, f, \mathscr{A}) \leq \mathrm{P}_n(T, f, \mathscr{A}) \leq e^{\lambda(n)\varepsilon} \mathrm{Q}_n(T, f, \mathscr{A}).$$

By taking the logarithm, dividing by $\lambda(n)$ and taking the limit for $n_\alpha$, we get equation (1). □

**Definition 3.6** (Topological Pressure). *For a G-action $T : X \to X$, and a one-point uniformly continuous function $f : X \to \mathbb{R}$, define $\mathrm{P}(T, f)$ to be the quantity presented in lemma 3.5. That is,*

$$\begin{aligned}
\mathrm{P}(T, f) &= \lim_{\mathscr{A}:\ admissible\ cover} \mathrm{P}(T, f, \mathscr{A}) \\
&= \lim_{\mathscr{A}:\ admissible\ cover} \mathrm{Q}(T, f, \mathscr{A}) \\
&= \sup_{\mathscr{A}:\ admissible\ cover} \mathrm{Q}(T, f, \mathscr{A}).
\end{aligned}$$

## 3.2 Using Separated and Spanning Sets

Traditionally, concepts of *separated sets* and *spanning sets* are defined using a metric or, more generaly, a uniform structure. We shall use *admissible covers.* Given a cover $\mathscr{A}$, let's say that two points $x, y \in X$ are $\mathscr{A}$-close if there is a $B \in \mathscr{A}$ such that $x, y \in B$.

A subset $E \subset X$ is said to be $\mathscr{A}$-spanning if given $x \in X$ there is $y \in E$ such that $x$ is $\mathscr{A}$-close to $y$. While $E$ is said to be a $\mathscr{A}$-separated if no distinct elements of $E$ are $\mathscr{A}$-close. It is easily shown that every maximal $\mathscr{A}$-separated set is $\mathscr{A}$-spanning. In fact, if $E$ is $\mathscr{A}$-separated but not $\mathscr{A}$-spanning, then there is an $x \in X$ that is not $\mathscr{A}$-close to any point of $E$. Therefore, $E \cup \{x\}$ is still $\mathscr{A}$-separated.

**Definition 3.7.** *Let $T : X \to X$ be a G-action. Given a cover $\mathscr{A}$, define*

$$\mathrm{G}_n(T, f, \mathscr{A}) = \inf \left\{ \sum_{x \in E} \mathrm{e}^{f_n(x)} \,\middle|\, E \text{ is } \mathscr{A}^n\text{-spanning} \right\}$$

$$\mathrm{S}_n(T, f, \mathscr{A}) = \sup \left\{ \sum_{x \in E} \mathrm{e}^{f_n(x)} \,\middle|\, E \text{ is } \mathscr{A}^n\text{-separated} \right\}.$$

*and*

$$\begin{aligned} \mathrm{G}(T, f, \mathscr{A}) &= \lim_{\alpha} \frac{1}{\lambda(n_\alpha)} \log \mathrm{G}_{n_\alpha}(T, f, \mathscr{A}) \\ \mathrm{S}(T, f, \mathscr{A}) &= \lim_{\alpha} \frac{1}{\lambda(n_\alpha)} \log \mathrm{S}_{n_\alpha}(T, f, \mathscr{A}). \end{aligned}$$

It is clear that

$$\mathscr{A} \prec \mathscr{B} \Rightarrow \begin{cases} \mathrm{G}(T, f, \mathscr{A}) \leq \mathrm{G}(T, f, \mathscr{B}) \\ \mathrm{S}(T, f, \mathscr{A}) \leq \mathrm{S}(T, f, \mathscr{B}). \end{cases}$$

Therefore, the limits in the following definition exist and are equal to the supremum.

**Definition 3.8.** *Let $T : X \to X$ a G-action. Define*

$$\begin{aligned} \mathrm{G}(T, f) &= \lim_{\mathscr{A}:\ \textit{admissible cover}} \mathrm{G}(T, f, \mathscr{A}) \\ \mathrm{S}(T, f) &= \lim_{\mathscr{A}:\ \textit{admissible cover}} \mathrm{S}(T, f, \mathscr{A}). \end{aligned}$$

## 3.3 Properties

**Lemma 3.9.** *Let $T : X \to X$ be a G-action over a locally compact Hausdorff topological space $X$. And let $f : X \to \mathbb{R}$ be one-point uniformly continuous. Then,*

$$\mathrm{S}(T, f) \leq \mathrm{P}(T, f).$$

*Proof.* Given an admissible cover $\mathscr{A}$, let $E \subset X$ be $\mathscr{A}^n$-separated. Given a subcover $\mathscr{B} \subset \mathscr{A}^n$, for each $x \in X$, take $B_x \in \mathscr{B}$ with $x \in B_x$. Since $B_x \neq B_y$ for any distinct pair $x, y \in E$,

$$\sum_{x \in E} e^{f_n(x)} \leq \sum_{x \in E} \sup e^{f_n(B_x)} \leq \sum_{B \in \mathscr{B}} \sup e^{f_n(B)}.$$

By taking the infimum over $\mathscr{B} \subset \mathscr{A}^n$,

$$\sum_{x \in E} e^{f_n(x)} \leq \mathrm{P}_n(T, f, \mathscr{A}).$$

Taking the supremum for the $\mathscr{A}^n$-separated sets $E$,

$$\mathrm{S}_n(T, f, \mathscr{A}) \leq \mathrm{P}_n(T, f, \mathscr{A}).$$

Taking the limit for $n = n_\alpha$,

$$\mathrm{S}(T, f, \mathscr{A}) \leq \mathrm{P}(T, f, \mathscr{A}) \leq \mathrm{P}(T, f).$$

Since $\mathscr{A}$ was arbitrary,

$$\mathrm{S}(T, f) \leq \mathrm{P}(T, f).$$

□

**Lemma 3.10.** *Let $T : X \to X$ be a $G$-action over a locally compact Hausdorff topological space $X$. And let $f : X \to \mathbb{R}$ be one-point uniformly continuous. Then,*

$$\mathrm{P}(T, f) \le \mathrm{G}(T, f).$$

*Proof.* Let $\mathscr{A}$ be an admissible cover. Let $E$ be $\mathscr{A}^n$-spanning. For each $x \in E$, pick a $B_x \in \mathscr{A}^n$ such that $x \in B_x$. The fact that $E$ is $\mathscr{A}^n$-spanning implies that $\mathscr{B} = \left\{B_x \,\middle|\, x \in E\right\}$ is a cover. Therefore,

$$\begin{aligned} \mathrm{Q}_n(T, f, \mathscr{A}) &\le \sum_{x \in E} \inf e^{f_n(B_x)} \\ &\le \sum_{x \in E} e^{f_n(x)}. \end{aligned}$$

On may now take the infimum over the $\mathscr{A}^n$-generating sets $E$ to find

$$\mathrm{Q}_n(T, f, \mathscr{A}) \le \mathrm{G}_n(T, f, \mathscr{A}).$$

Now, just take the limit for $n = n_\alpha$ and then for $\mathscr{A}$. □

**Proposition 3.11.** *For a $G$-action $T : X \to X$ over a locally compact Hausdorff space and a one-point uniformly continuous function $f : X \to \mathbb{R}$,*

$$\mathrm{P}(T, f) = \mathrm{G}(T, f) = \mathrm{S}(T, f).$$

*Moreover, none of those quantities depend on the choice of the universal subnet $n_\alpha$.*

*Proof.* Since a maximal $\mathscr{A}^n$-separated set is $\mathscr{A}^n$-generating,

$$\mathrm{G}(T, f) \le \mathrm{S}(T, f).$$

With lemmas 3.9 and 3.10,

$$\mathrm{P}(T, f) \le \mathrm{G}(T, f) \le \mathrm{S}(T, f) \le \mathrm{P}(T, f).$$

Since $\mathrm{P}(T, f)$ does not depend on the choice of $n_\alpha$, none of them do. □

It shall be more convenient to work with $f \in C_0(X)$, instead of an arbitrary one-point continuous $f : X \to \mathbb{R}$. If $f$ is one-point continuous, then $f - f(\infty) \in C_0(X)$.

**Lemma 3.12.** *Let $T : X \to X$ be a Hausdorff and locally compact dynamical system. Let $f : X \to \mathbb{R}$ be one-point uniformly continuous and $c \in \mathbb{R}$. Then,*

$$\mathrm{P}(T, f + c) = \mathrm{P}(T, f) + c.$$

*Proof.* The lemma is evident from the fact that, for $E \subset X$,

$$\log \sum_{x \in E} \mathrm{e}^{f(x) + c} = c + \log \sum_{x \in E} \mathrm{e}^{f(x)} .$$

□

# 4 Variational Principle

Inspired by what has been done for the compact case, we demonstrate a *variational principle* for the pressure of a topological system $T : X \to X$, where $X$ is not assumed to be compact but it is just assumed to have a one-point compactification $X^*$. This does not imply that $T$ can be itself extended to a topological dynamical system over $X^*$.

We use the preparations made in sections 2 and 3 in order to adapt Misiurewicz's demonstration of the variational principle. Misiurewicz's original article is [6]. We shall follow the more didatic presentation of the variational principle presented in [8, sections 10.3 and 10.4]. A similar presentation can also be found in [10, chapter 9].

We shall deal with the supremum of $P_\mu(T, f)$ over all $T$-invariant Radon probability measures over $X$. However, there might happen that no such a probability measure exists.

**Theorem 4.1.** *Let $T : X \to X$ be a topological dynamical system, and let $f : X \to \mathbb{R}$ be* one-point uniformly continuous*. Then,*

$$\mathrm{P}(T, f) = \sup_\mu P_\mu(T, f),$$

*where the supremum is taken over all $T$-invariant Radon probability measures. If there is no $T$-invariant Radon probability measure,*

$$\mathrm{P}(T, f) = f(\infty).$$

*Proof.* If we assume that the theorem holds for $f \in C_0(X)$, then lemmas 2.23 and 3.12 show that

$$\begin{aligned}\mathrm{P}(T, f) &= \mathrm{P}(T, f - f(\infty)) + f(\infty) \\ &= \sup_\mu P_\mu(T, f - f(\infty)) + f(\infty)\end{aligned}$$

When there are no $T$-invariant probabilities, the last term in this equation is just $f(\infty)$. When there are $T$-invariant probabilities, according to lemma 2.23,

$$\begin{aligned}\sup_\mu P_\mu(T, f - f(\infty)) + f(\infty) &= \sup_\mu P_\mu(T, f) - f(\infty) + f(\infty) \\ &= \sup_\mu P_\mu(T, f).\end{aligned}$$

Therefore, we shall assume $f \in C_0(X)$.

The theorem will be demonstrated if we show that:

1. For any $T$-invariant Radon measure $\mu$ with $0 \leq \mu(X) \leq 1$,

   $$P_\mu(T, f) \leq \mathrm{P}(T, f).$$

2. Given an *admissible cover* $\mathscr{A}$, there is a $T$-invariant Radon measure $\mu$, with $0 \leq \mu(X) \leq 1$, such that

   $$\mathrm{S}(T, f, \mathscr{A}) \leq P_\mu(T, f).$$

These claims are the contents of the following two subsections. □

## 4.1 Topological Pressure is an Upper Bound.

**Proposition 4.2.** *Let $T : X \to X$ be a dynamical system. Let $f : X \to \mathbb{R}$ be* one-point uniformly continuous*. And let $\mu$ be a $T$-invariant measure with $0 \leq \mu(X) \leq 1$. Then,*

$$P_\mu(T, f) \leq \mathrm{P}(T, f).$$

*Proof.* The claim is evident if $\mu = 0$. By lemma 2.24, we may assume $\mu$ is a probability measure. We shall show that for any $n \in G$,

$$P_\mu(T^n, f_n) \leq \lambda(n)\mathrm{P}(T, f) + 2 + \log 2. \tag{2}$$

And then, lemma 2.30 implies that

$$\begin{aligned} P_\mu(T, f) &= \frac{1}{\lambda(n)} P_\mu(T^n, f_n) \\ &\leq \mathrm{P}(T, f) + \frac{2 + \log 2}{\lambda(n)} \to \mathrm{P}(T, f), \end{aligned}$$

concluding the demonstration.

Notice that $f_n$ might not be one-point uniformly continuous, and therefore, we do not talk about $\mathrm{P}(T^n, f_n)$. From now on, we fix $n$ and show the validity of inequation (2).

According to proposition 2.26, we have to show that given an admissible partition $\mathscr{K}$,

$$P_\mu(T^n, f_n, \mathscr{K}) < \lambda(n)\mathrm{P}(T, f) + 2 + \log 2.$$

So, we just have to prove that there is an admissible cover $\mathscr{B}$ such that

$$P_\mu(T^n, f_n, \mathscr{K}) \leq \lambda(n)\mathrm{Q}(T, f, \mathscr{B}) + 2 + \log 2. \tag{3}$$

Let $\mathscr{A}$ be the strongly admissible cover associated to $\mathscr{K}$, from lemma 2.13. Let $\varepsilon = 1/\lambda(n)$, and consider $\mathscr{B}_{f,\varepsilon}$ from corollary 2.16. The cover

$$\mathscr{B} = \mathscr{A} \vee \mathscr{B}_{f,\varepsilon},$$

as an open refinement of the *admissible cover* $\mathscr{B}_{f,\varepsilon}$, is *admissible.* So, $\mathscr{B}$ is an *admissible cover* that refines $\mathscr{A}$ and is such that

$$x, y \in B \in \mathscr{B}^{mn} \Rightarrow f_{mn}(x) \leq f_{mn}(y) + \frac{\lambda(mn)}{\lambda(n)} = f_{mn}(y) + \lambda(m). \tag{4}$$

With $\mathscr{B}$ properly chosen, we demonstrate the validity of inequality (3). And since we are working with $T$ and $T^n$ at the same time, let's agree that whenever the transformation is omitted, it is assumed to be $T$.

Given $m \in G$, let $\mathscr{D} \subset \mathscr{B}^{mn}$ be any subcover. And notice that

$$\begin{aligned} (\mathscr{K}^n)^m_{T^n} &\prec \mathscr{K}^{mn} \\ (f_n)_{T^n,m} &= f_{mn}. \end{aligned}$$

Given $C \in \mathscr{K}^m_{T^n}$, let $x_C \in C$ be such that

$$\sup f_{mn}(C) \leq f_{mn}(x_C) + 1,$$

Also, and choose $B_C \in \mathscr{D}$ such that $x_C \in B_C$. Then, equation (4) implies that

$$\begin{aligned}\sup f_{mn}(C) &\leq f_{mn}(x_C) + 1 \\ &\leq \inf f_{mn}(B_C) + \lambda(m) + 1.\end{aligned}$$

For each $B \in \mathscr{D}$, let

$$c_B = \#\left\{C \in \mathscr{K}_{T^n}^m \,\middle|\, B_C = B\right\}.$$

Since $\mathscr{A} \prec \mathscr{B} \prec \mathscr{D}$, lemma 2.14 implies that

$$c_B \leq 2^{\lambda(m)}.$$

Using lemma 2.31 with $T^n$ instead of $T$, $f_n$ instead of $f$, and $m$ instead of $n$,

$$\begin{aligned}\int f_n \,\mathrm{d}\mu + \frac{1}{\lambda(m)} H_\mu(\mathscr{K}_{T^n}^m) &\leq \frac{1}{\lambda(m)} \log \sum_{C \in \mathscr{K}_{T^n}^m} \mathrm{e}^{\sup f_{mn}(C)} \\ &\leq \frac{1}{\lambda(m)} \log(\mathrm{e}^{\lambda(m)+1} \sum_{C \in \mathscr{K}_{T^n}^m} \mathrm{e}^{\inf f_{mn}(B_C)}) \\ &= \frac{\lambda(m)+1}{\lambda(m)} + \frac{1}{\lambda(m)} \log \sum_{C \in \mathscr{K}_{T^n}^m} \mathrm{e}^{\inf f_{mn}(B_C)} \\ &\leq \frac{\lambda(m)+1}{\lambda(m)} + \frac{1}{\lambda(m)} \log \sum_{B \in \mathscr{D}} c_B \,\mathrm{e}^{\inf f_{mn}(B)} \\ &\leq \frac{\lambda(m)+1}{\lambda(m)} + \frac{1}{\lambda(m)} \log(2^{\lambda(m)} \sum_{B \in \mathscr{D}} \mathrm{e}^{\inf f_{mn}(B)}) \\ &= \frac{\lambda(m)+1}{\lambda(m)} + \log 2 + \frac{\lambda(n)}{\lambda(mn)} \log \sum_{B \in \mathscr{D}} \mathrm{e}^{\inf f_{mn}(B)}.\end{aligned}$$

Taking the infimum for every subcover $\mathscr{D} \subset \mathscr{B}^{mn}$,

$$\begin{aligned}&\int f_n \,\mathrm{d}\mu + \frac{1}{\lambda(m)} H_\mu(\mathscr{K}_{T^n}^m) \\ &\qquad \leq \frac{\lambda(m)+1}{\lambda(m)} + \log 2 + \frac{\lambda(n)}{\lambda(mn)} \log \mathrm{Q}_{nm}(T, f, \mathscr{B}).\end{aligned}$$

Now, take the limit on $m$ to get

$$\begin{aligned}\int f_n \,\mathrm{d}\mu + h_\mu(T^n), \mathscr{K} &\leq 1 + \log 2 + \lim_{m \in G} \frac{\lambda(n)}{\lambda(nm)} \log \mathrm{Q}_{nm}(T, f, \mathscr{B}) \\ &= 1 + \log 2 + \lambda(n)\mathrm{Q}(T, f, \mathscr{B}) \\ &\leq 1 + \log 2 + \lambda(n)\mathrm{P}(T, f).\end{aligned}$$

This is inequality (3), and thus we conclude the proof. □

## 4.2 Topological Pressure is a Lower Bound

In this subsection, we prove the following proposition, which is an adaption of what is done in [8, subsection 10.4.2], using techniques similar to those used in [3, theorem 3.1].

**Proposition 4.3.** *Let $T : X \to X$ be a Hausdorff locally compact dynamical system. Then, for any $f \in C_0(X)$ and any* admissible cover $\mathscr{A}$, *there exists a $T$-invariant Radon measure $\mu$, with $0 \leq \mu(X) \leq 1$, such that*

$$\mathrm{S}(T, f, \mathscr{A}) \leq P_\mu(T, f).$$

First, use lemma 2.7 to get a compact Hausdorff extension $S : Z \to Z$ for $T$. According to lemma 2.25, the demonstration will be complete if we find an $S$-invariant probability measure $\mu$ over $Z$ and a $\pi$-measurable partition $\mathscr{Z}$, such that

$$\mathrm{S}(T, f, \mathscr{A}) \leq P_\mu(S, g, \mathscr{Z}). \tag{5}$$

where $g = f \circ \pi$.

### 4.2.1 The construction of $\mu$

Consider an *admissible cover* $\mathscr{A}$. For each $n \in G$, let $E_n \subset X$ be $\mathscr{A}^n$-separated and such that

$$\frac{1}{2}\mathrm{S}_n(T, f, \mathscr{A}) \leq \sum_{x \in E_n} \mathrm{e}^{f_n(x)}.$$

Call the right hand side quantity $S_n$. So,

$$\frac{1}{2}\mathrm{S}_n(T, f, \mathscr{A}) \leq S_n.$$

Then, define over $Z$ the probability measure

$$\sigma_n = \frac{1}{S_n} \sum_{x \in E_n} \mathrm{e}^{g_n(x)} \, \delta_x,$$

where $\delta_x$ is the Dirac measure with support in $x$. And then, define

$$\mu_n = \frac{1}{\lambda(n)} \sum_{k \in \Lambda(n)} \sigma_n \circ S^{-k}.$$

Since the set of probability measures over $Z$ is compact in the weak-$*$ topology (see subsection 2.1), any universal subnet is convergent to a probability measure. To avoid introducing a new variable, let's just use our already chosen $n_\alpha$. Then, $\mu_{n_\alpha}$ converges to a certain probability measure $\mu$. In special, by the *Portmanteau Theorem* (proposition 2.2), if $C$ is measurable with $\mu(\partial C) = 0$,

$$\mu_{n_\alpha}(C) \to \mu(C).$$

We claim that $\mu$ is $S$-invariant. In fact, for any continuous $\varphi : Z \to \mathbb{R}$ and any $m \in G$, $\varphi \circ S^m$ is also continuous. Therefore,

$$\begin{aligned} \int \varphi \, \mathrm{d}(\mu_{n_\alpha} \circ S^{-m}) &= \int \varphi \circ S^m \, \mathrm{d}\mu_{n_\alpha} \\ &\to \int \varphi \circ S^m \, \mathrm{d}\mu \\ &= \int \varphi \, \mathrm{d}(\mu \circ S^{-m}). \end{aligned}$$

That is, $\mu_{n_\alpha} \circ S^{-m} \to \mu \circ S^{-m}$. On the other hand, let $\Gamma_n = \Lambda(n) \,\Delta\, (m + \Lambda(n))$ and use lemma 2.3 to get

$$\begin{aligned} |\int \varphi \,\mathrm{d}(\mu_{n_\alpha} - \mu_{n_\alpha} \circ S^{-m})| &= \frac{1}{\lambda(n_\alpha)} |\sum_{k \in \Lambda(n_\alpha)} \int \varphi \,\mathrm{d}(\sigma_{n_\alpha} \circ S^{-k} - \sigma_{n_\alpha} \circ S^{-(m+k)})| \\ &= \frac{1}{\lambda(n_\alpha)} \sum_{k \in \Gamma_{n_\alpha}} |\int \varphi \,\mathrm{d}\sigma_{n_\alpha} \circ S^{-k}| \\ &\leq \frac{\#\Gamma_{n_\alpha}}{\lambda(n_\alpha)} \|\varphi\|_\infty \to 0. \end{aligned}$$

That is, $\lim \mu_{n_\alpha} = \lim \mu_{n_\alpha} \circ S^{-m}$. Therefore,

$$\mu = \lim \mu_{n_\alpha} = \lim \mu_{n_\alpha} \circ S^{-m} = \mu \circ S^{-m}.$$

**Lemma 4.4.** $\frac{1}{\lambda(n)} \int g_n \,\mathrm{d}\sigma_n = \int g \,\mathrm{d}\mu_n$.

*Proof.* In fact,

$$\begin{aligned} \frac{1}{\lambda(n)} \int g_n \,\mathrm{d}\sigma_n &= \frac{1}{\lambda(n)} \sum_{k \in \Lambda(n)} \int g \circ S^k \,\mathrm{d}\sigma_n \\ &= \frac{1}{\lambda(n)} \sum_{k \in \Lambda(n)} \int g \,\mathrm{d}\sigma_n \circ S^{-k} \\ &= \int g \,\mathrm{d}\mu_n. \end{aligned}$$

□

### 4.2.2 The construction of $\mathscr{Z}$

We construct a suitable $\pi$-measurable partition $\mathscr{Z}$, in such a way that inequation (5) holds. Lemma 2.15 says that $\mathscr{A}$ is the restriction to $X$ of some open cover of $X^*$, $\tilde{\mathscr{A}}$. Use lemma 2.17 with $\mu \circ \pi^{-1}$ to get a finite open refinement $\mathscr{B}$ for $\tilde{\mathscr{A}}$ such that every member of $\mathscr{B}$ has a $(\mu \circ \pi^{-1})$-null border. Since $X^*$ is compact, we can assume that $\mathscr{B} = \{B_1, \dots, B_l\}$ is finite and has no proper subcover. Then,

$$Z_j = \pi^{-1}(B_j) \setminus (\pi^{-1}(B_1) \cup \dots \cup \pi^{-1}(B_{j-1}))$$

form a $\pi$-measurable partition of $Z$, $\mathscr{Z} = \{Z_1, \dots, Z_l\}$.

**Lemma 4.5.** *For each $C \in \mathscr{Z}^n$,*

$$\mu_{n_\alpha}(C) \to \mu(C).$$

*Proof.* Since $S^m$ is continuous for every $m \in G$, the border operator $\partial$ possesses the following properties.

1. $\partial A = \partial(A^c)$.
2. $\partial(A_1 \cap \dots \cap A_r) \subset \partial A_1 \cup \dots \cup \partial A_r$.
3. $\partial S^{-m}(A) \subset S^{-m}(\partial A)$.

From items (1) and (2), each $Z_j = B_j \cap {B_1}^c \cap \cdots \cap {B_{j-1}}^c$ in $\mathscr{X}$ has border with null measure. And from items (2) and (3), the same is true for any set $C \in \mathscr{X}^n$. Therefore, *Portmanteau Theorem* (proposition 2.2) implies that for those sets, $\mu_{n_\alpha}(C) \to \mu(C)$. □

We have constructed a partition $\mathscr{X}$ that can be used with lemma 2.25. And also, since $E_n$ is $\mathscr{A}^n$-separated, $E_n$ is $\mathscr{X}^n$-separated as well.

### 4.2.3 Comparing with separated sets

Having constructed $\mu$ and $\mathscr{X}$, it remains to show that inequation (5) holds.

**Lemma 4.6.** $\log S_n = H_{\sigma_n}(\mathscr{X}^n) + \int g_n \, \mathrm{d}\sigma_n$.

*Proof.* By construction, since $\tilde{\mathscr{A}} \prec \mathscr{X}$, each $C \in \mathscr{X}^n$ contains at most one element $x \in E_n$. That is, either $\sigma_n(C) = 0$ or $\sigma_n(C) = \frac{\mathrm{e}^{g_n(x)}}{S_n}$. Therefore,

$$\begin{aligned}
H_{\sigma_n}(\mathscr{X}^n) + \int g_n \, \mathrm{d}\sigma_n &= \sum_{x \in E_n} \sigma_n(\{x\})(\log \frac{1}{\sigma_n(\{x\})} + g_n(x)) \\
&= \sum_{x \in E_n} \frac{\mathrm{e}^{g_n(x)}}{S_n} \log \frac{\mathrm{e}^{g_n(x)}}{\mathrm{e}^{g_n(x)}/S_n} \\
&= \sum_{x \in E_n} \frac{\mathrm{e}^{g_n(x)}}{S_n} \log S_n \\
&= \log S_n.
\end{aligned}$$

□

Passing from $\sigma_n$ to $\mu_n$ follows the same steps as in the compact case. Notice that for any measurable finite partition $\mathscr{D}$, [3, lemma 2.7] implies that

$$\sum_{k \in \Lambda(n)} \frac{1}{\lambda(n)} H_{\sigma_n \circ S^{-k}}(\mathscr{D}) \leq H_{\mu_n}(\mathscr{D}). \tag{6}$$

And, fixed $n, q \in G$ with $q < n$, then, for each $k \in \Lambda(q)$,

$$\mathscr{X}^n = \mathscr{X}^{\Gamma_n(q,k)} \vee \Big( \bigvee_{p \in \Lambda_n(q,k)} \mathscr{X}^{p+\Lambda(q)} \Big)$$

So, for every $k \in \Lambda(q)$ we can use lemma 2.19 to obtain

$$\begin{aligned}
H_{\sigma_n}(\mathscr{X}^n) &\leq H_{\sigma_n}(\mathscr{X}^{\Gamma_n(q,k)}) + \sum_{p \in \Lambda_n(q,k)} H_{\sigma_n}(\mathscr{X}^{p+\Lambda(q)}) \\
&= H_{\sigma_n}(\mathscr{X}^{\Gamma_n(q,k)}) + \sum_{p \in \Lambda_n(q,k)} H_{\sigma_n \circ S^{-p}}(\mathscr{X}^q) \\
&\leq \#\Gamma_n(q,k) \log \#\mathscr{X} + \sum_{p \in \Lambda_n(q,k)} H_{\sigma_n \circ S^{-p}}(\mathscr{X}^q).
\end{aligned}$$

Summing for all elements $k \in \Lambda(q)$,

$$\begin{aligned}
\lambda(q) H_{\sigma_n}(\mathscr{X}^n) &\leq (\sum_{k\in\Lambda(q)} \#\Gamma_n(q,k)) \log \#\mathscr{X} + \\
&\quad + \sum_{k\in\Lambda(q)} (\sum_{p\in\Lambda_n(q,k)} H_{\sigma_n \circ S^{-p}}(\mathscr{X}^q)) \\
&\leq (\sum_{k\in\Lambda(q)} \#\Gamma_n(q,k)) \log \#\mathscr{X} + \sum_{p\in\Lambda(n)} H_{\sigma_n \circ S^{-p}}(\mathscr{X}^q) \\
&\leq (\sum_{k\in\Lambda(q)} \#\Gamma_n(q,k)) \log \#\mathscr{X} + \lambda(n) H_{\mu_n}(\mathscr{X}^q),
\end{aligned} \tag{7}$$

where the last inequality comes from (6). Since $\mu_{n_\alpha}(C) \to \mu(C)$ for any $C \in \mathscr{X}^q$,

$$H_{\mu_{n_\alpha}}(\mathscr{X}^q) \to H_\mu(\mathscr{X}^q).$$

Therefore, using lemma 4.6, followed by inequation (7) and lemma 4.4; and then, lemma 2.3,

$$\begin{aligned}
\mathrm{S}(T, f, \mathscr{A}) &= \lim_\alpha \frac{1}{\lambda(n_\alpha)} \log \mathrm{S}_{n_\alpha}(T, f, \mathscr{A}) \\
&= \lim_\alpha \frac{1}{\lambda(n_\alpha)} \log \frac{1}{2} \mathrm{S}_{n_\alpha}(T, f, \mathscr{A}) \\
&\leq \lim_\alpha \frac{1}{\lambda(n_\alpha)} \log S_{n_\alpha} \\
&= \lim_\alpha \frac{1}{\lambda(n_\alpha)} H_{\sigma_{n_\alpha}}(\mathscr{X}^{n_\alpha}) + \lim_\alpha \frac{1}{\lambda(n_\alpha)} \int g_{n_\alpha} \, \mathrm{d}\sigma_{n_\alpha} \\
&= \frac{1}{\lambda(q)} \lim_\alpha \frac{\lambda(q) H_{\sigma_{n_\alpha}}(\mathscr{X}^{n_\alpha})}{\lambda(n_\alpha)} + \lim_\alpha \int g \, \mathrm{d}\mu_{n_\alpha} \\
&\leq \frac{1}{\lambda(q)} \lim_\alpha (\sum_{k\in\Lambda(q)} \frac{\#\Gamma_{n_\alpha}(q,k)}{\lambda(n_\alpha)}) \log \#\mathscr{X} + \\
&\quad + \frac{1}{\lambda(q)} \lim_\alpha H_{\mu_{n_\alpha}}(\mathscr{X}^q) + \lim_\alpha \int g \, \mathrm{d}\mu_{n_\alpha} \\
&= 0 + \frac{1}{\lambda(q)} H_\mu(\mathscr{X}^q) + \int g \, \mathrm{d}\mu.
\end{aligned}$$

Now, taking the limit for $q$,

$$\mathrm{S}(T, f, \mathscr{A}) \leq P_\mu(S, g, \mathscr{X}),$$

concluding the demonstration of inequation (5), and therefore, concluding the demonstration of proposition 4.3.

# 5 Entropy and Pressure Leakage

Different attempts can be made in order to define *topological pressure* for non compact systems. In this section, we shall illustrate an example that shows how some alternatives would fail to satisfy the *variational principle*. One possible interpretation is that those alternatives detect the entropy (complexity) present in the border of the compactified version of the system.

In the example that follows, we show the importance of the concept of *admissible cover* when defining *topological pressure*. If instead of *admissible covers*, we use *finite open covers* the variational principle does not hold. And if, instead of counting $\mathscr{A}^n$-separated points, we count $(n, \varepsilon)$-separated points for the unique uniform structure compatible with the compactified version of the system, the *variational principle* does not hold, either.

## 5.1 Dynamical System over the Open Disk

First, we shall define the dynamical system of our example. Let

$$\begin{aligned} X &= \left\{ z \in \mathbb{C} \,\middle|\, |z| < 1 \right\} \\ Z &= \left\{ z \in \mathbb{C} \,\middle|\, |z| \leq 1 \right\} \\ \partial X &= \left\{ z \in \mathbb{C} \,\middle|\, |z| = 1 \right\} \end{aligned}$$

be the open and closed unit disks in $\mathbb{C}$ and their border. Consider the dynamical system

$$\begin{array}{rccc} S: & Z & \to & Z \\ & r\,\mathrm{e}^{i\theta} & \mapsto & r^{\frac{r+1}{2}}\,\mathrm{e}^{i2\theta} \end{array}$$

and consider its restriction to $X$, $T : X \to X$, and to $\partial X$,

$$\begin{array}{rccc} K: & \partial X & \to & \partial X \\ & \mathrm{e}^{i\theta} & \mapsto & \mathrm{e}^{i2\theta} \end{array}.$$

It is a simple fact that

$$h(T) = 0.$$

In fact, Poincaré's Recurrence Theorem [8, theorem 1.2.4] implies that the only invarian probability measure for $T$ is Dirac's $\delta_0$. And this implies, using the *Variational Principle* (theorem 4.1), that $h(T) = 0$.

On the other hand, the dynamical system $K$, over $\partial X$, is such that

$$h(K) = \log 2.$$

## 5.2 Counting Separated Points

In this subsection, we shall consider $(n, \varepsilon)$-separated sets, like defined in [6]. Continuing with our example, take the usual euclidean metric $d$ over $Z$. For $\varepsilon_n \to 0$, there is a sequence of $(n, \varepsilon_n)$-separated sets $E_n$, such that

$$\lim_{n\to\infty} \frac{1}{n} \log \# E_n \to \log 2.$$

Now, if we take points sufficiently close to those of $E_n$, e can produce $F_n \subset X$ that is still $(n, \varepsilon_n)$-separated. That is, with the euclidean metric, separated sets would give us a non-zero entropy for $T$, because

$$\lim_{n\to\infty} \frac{1}{n} \log \# F_n \to \log 2.$$

Using $\mathscr{A}^n$ separated sets for an *admissible cover* $\mathscr{A}$ is just like using a pseudo-metric over $Z$ that gives diameter 0 to $\partial X$. In this case, $E_n$ would never be $(n, \varepsilon_n)$-separated.

### 5.3 Using Finite Covers

An argument similar to that used in subsection 5.2 works for the idea of using finite covers instead of *admissible covers.* There is a sequence of minimal open covers $\mathscr{B}_n$ such that

$$\lim_{n\to\infty} \frac{1}{n} \log N(\mathscr{B}_n^n) \to \log 2.$$

To each set $B \subset \partial X$, define the *"pizza slice"*

$$A_B = \left\{rx \,\middle|\, x \in B,\, 0 < r < 1\right\},$$

and let $D = \left\{x \in X \,\middle|\, |x| < 1/2\right\}$. Then,

$$\mathscr{A}_n = \{D\} \cup \left\{A_B \,\middle|\, B \in \mathscr{B}_n\right\}$$

is a finite open cover for $X$ such that

$$N(\mathscr{A}_n^n) = N(\mathscr{B}_n^n) + 1.$$

And again, this implies that entropy defined using finite covers would not satisfy the *variational principle*, as it would assign a non-zero entorpy to $T$.